\documentclass[12pt, a4paper]{amsart}
\pdfoutput=1
\usepackage{amsmath,amsfonts,amssymb,amsthm,color}  
\usepackage{mathrsfs}

\usepackage[utf8]{inputenc}
\usepackage{graphicx}

\theoremstyle{plain}
\newtheorem{theorem}{Theorem}[section]

\newtheorem{lemma}[theorem]{Lemma}
\newtheorem{proposition}[theorem]{Proposition}

\theoremstyle{definition}

\newtheorem{remark}[theorem]{Remark}

\numberwithin{equation}{section}
\newtheorem*{theorem*}{Theorem}

\def\Xint#1{\mathchoice
{\XXint\displaystyle\textstyle{#1}}
{\XXint\textstyle\scriptstyle{#1}}
{\XXint\scriptstyle\scriptscriptstyle{#1}}
{\XXint\scriptscriptstyle\scriptscriptstyle{#1}}
\!\int}
\def\XXint#1#2#3{{\setbox0=\hbox{$#1{#2#3}{\int}$}
\vcenter{\hbox{$#2#3$}}\kern-.5\wd0}}

\def\dashint{\Xint-}

\newcommand{\dy}{\, \mathrm{d}y}
\newcommand{\dx}{\, \mathrm{d}x}
\newcommand{\dt}{\, \mathrm{d}t}

\DeclareMathOperator{\dive}{div}

\DeclareMathOperator{\BMO}{BMO}
\DeclareMathOperator{\PBMO}{PBMO}

\providecommand{\abs}[1]{ \lvert#1  \rvert}
\providecommand{\norm}[1]{ \lVert#1  \rVert}

\providecommand{\s}[1]{\mathscr{#1}}

\title[Parabolic BMO and the maximal operator]{Parabolic BMO and the forward-in-time maximal operator}
\author{Olli Saari}
\begin{document}
\begin{abstract}
We study if the parabolic forward-in-time maximal operator is bounded on parabolic $\BMO$. It turns out that for non-negative functions the answer is positive, but the behaviour of sign changing functions is more delicate. The class parabolic $\BMO$ and the forward-in-time maximal operator originate from the regularity theory of nonlinear parabolic partial differential equations. In addition to that context, we also study the question in dimension one.
\end{abstract}

\address{Olli Saari,	
	Department of Mathematics and Systems Analysis, 
	Aalto University School of Science,
	FI-00076 Aalto, 
	Finland
} \email{olli.saari@aalto.fi} 



\thanks{The author is supported by the V\"ais\"al\"a Foundation.}

\maketitle


\section{Introduction}
The Hardy-Littlewood maximal operator maps functions of bounded mean oscillation back to BMO. This is a classical result of Bennett, DeVore, and Sharpley \cite{BdVS1981}. In addition to the original approach, which is direct, some alternative proofs relying on the properties of the Muckenhoupt weights are available in \cite{CF1987,CN1995}. 

The present paper is devoted to studying the counterpart of the boundedness $M: \BMO \to \BMO$ in a context that comes from the regularity theory of parabolic partial differential equations \cite{FG1985,KK2007,Moser1964}. The class parabolic BMO is defined through a condition measuring mean oscillation in a special way. As opposed to the ordinary BMO space, a natural time lag appears in connection with parabolic BMO. It causes several challenges, many of which have only been addressed very recently; see \cite{KS2014,Saari2014}. Roughly speaking, the positive and negative parts of the deviation from a constant only satisfy bounds in disjoint regions of the space time. For a function $u$ of space and time to be in $\PBMO^{+}$, it suffices that 
\[
\sup_{R} \inf _{a \in \mathbb{R}} \left( \dashint_{R^{-}(\frac{1}{2})} (a-u)^{+} + \dashint_{R^{+}(\frac{1}{2})}(u-a)^{+} \right) < \infty;
\]
see Section \ref{subse:bmo} for precise definitions. The condition above leads to many properties similar to those of the ordinary BMO, but it also allows for the possibility of arbitrarily fast growth in the negative time direction. Consequently, the difference between parabolic BMO and its classical counterpart is remarkable. 

Principal examples of partial differential equations with connection to parabolic BMO are the heat equation and its generalizations, most notably the doubly nonlinear equation 
\begin{equation}
\label{eq:equ_intro}
\dfrac{\partial(\abs{u}^{p-2}u)}{\partial t} - \dive (\abs{\nabla u}^{p-2} \nabla u) = 0, \quad 1<p<\infty.
\end{equation}
Logarithms of positive local solutions to \eqref{eq:equ_intro}, possibly with measure data, belong to parabolic BMO \cite{KK2007,Moser1964,Saari2014}. This gives important examples of $\PBMO^{+}$ functions. A similar relation is also known to hold for elliptic partial differential equations and the ordinary BMO. Consequently, the functions in parabolic BMO relate to those in ordinary BMO in a manner analogous to how supersolutions of the heat equation relate to those of the Laplace equation.

\begin{figure}
\includegraphics[scale=1]{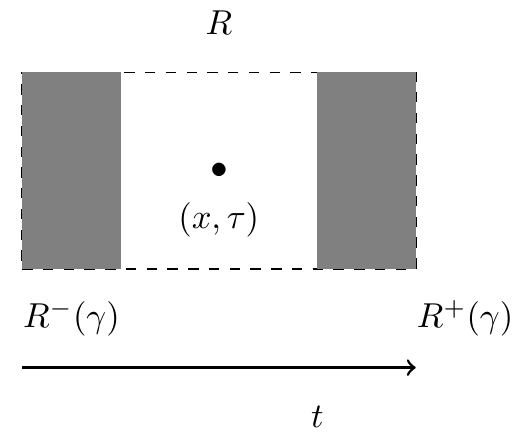}
\caption{\label{fig2:pararectangle}The sets $R^{\pm}(\gamma)$ in $\mathbb{R}^{n+1}$} 
\end{figure}

It was already noted by Moser \cite{Moser1964} and Trudinger \cite{Trudinger1968} in the 1960s that positive supersolutions of \eqref{eq:equ_intro} have their logarithms in a parabolic BMO class. This fact is well-known to people working in the regularity theory of partial differential equations; see \cite{Aimar1988,FG1985,KK2007}. However, the literature on parabolic $\BMO$ is still limited, and its history can be recounted quickly. The seminal papers of Moser \cite{Moser1964} and Trudinger \cite{Trudinger1968} established the connection between the partial differential equations and parabolic BMO. In addition, the parabolic John-Nirenberg inequality was proved there. The proof was simplified in \cite{FG1985} and extended to spaces of homogeneous type in \cite{Aimar1988}. Those papers date back to the 1980s. A method to derive global estimates from the local parabolic John-Nirenberg inequality was developed in \cite{Saari2014}, and it was also used to answer a question about summability of the supersolutions of \eqref{eq:equ_intro}. This was in 2014. More recent advances in the field are coupled with new trends in the theory of multidimensional one-sided weights; see \cite{FMRO2011,KS2014,KS2016_1,Ombrosi2005}. The techniques relevant in that context combine argumentation typical to the one-dimensional one-sided weight theory to that of harmonic and geometric analysis of partial differential equations. The usual challenge is to find a way to compensate the loss of many important tools such as the Besicovitch covering lemma and the standard maximal function techniques. 

Recall that the known results on parabolic $\BMO$ include a John--Nirenberg type inequality \cite{Aimar1988,Moser1964}, a Coifman-Rochberg type characterization through a special theory of weights \cite{KS2014}, and geometric local-to-global properties \cite{Saari2014}. The contribution of this paper is to show that the same operator of forward-in-time maximal averages (see Section \ref{sec:nota} and Figure \ref{fig2:pararectangle})
\[M^{\gamma+}f(x) = \sup_{\ell > 0} \dashint_{R^{+}(z,\ell,\gamma)} \abs{f}\] 
that generates the parabolic weight theory also maps the class of positive functions in parabolic $\BMO$ into itself. The boundedness of the parabolic forward-in-time maximal operator is consistent with the ``elliptic'' result in \cite{BdVS1981} when only positive functions are involved. Namely, the results coincide if we restrict our attention to functions with no time-dependency. In the time-dependent case, Theorem \ref{thm:main} is far more general. On the other hand, if we allow the functions under study to be sign-changing, a full analogue of the Bennett-DeVore-Sharpley result will be false in the parabolic context. Hence there are indigenously parabolic phenomena involved in our result. See Theorem \ref{thm:main} and the related discussion for more precise statements. 


We conclude the introduction by briefly describing the structure of the present paper.
It is divided into three main sections and an additional section discussing the one-dimensional analogue of the problem under study. Section \ref{sec:nota} is used to set up the notation and to introduce the operators and the function classes we study. In Section \ref{sec:preli}, we prove several auxiliary results such as a chain argument, and we also construct a special dyadic grid. These results are needed to solve problems that arise from the time-dependent nature of the main theorems. Once all the preparations have been carried out, we prove the main theorems, namely Lemma \ref{lemma:boundedness} asserting the boundedness of
\[M^{\gamma +} : \PBMO _{positive}^{+} \longrightarrow \PBMO_{positive}^{+}\] 
and Theorem \ref{thm:main} refining the result. At the end of the paper, we show how a simplified argument can be used the prove a slightly stronger result in dimension one. This is in the context of one-sided BMO spaces of Mart\'in-Reyes and de la Torre \cite{MRT1994}. To our best knowledge, also the one-dimensional result is new.

\vspace{0.5cm}

\noindent \textit{Acknowledgement.} The author would like to thank Juha Kinnunen for suggesting this problem. The author would also like to to thank Ioannis Parissis for many valuable comments on an earlier version of this manuscript.

\section{Notation and definitions}
\label{sec:nota}
We use standard notation. We mostly work in $\mathbb{R}^{n+1}$ with the last coordinate called time, the first ones space. The results also hold in space time cylinders that are sets of the form $\Omega \times \mathbb{R}$ with $\Omega \subset \mathbb{R}^{n}$ a bounded domain. The notation $A \lesssim B$ means that there is an uninteresting constant $C$ such that $A \leq CB$. We do not keep track of dependencies on dimension $n$, the growth type of the equation $p$ or parameters coming from the domain of definition $\Omega$. It is clear what $\gtrsim$ and $\eqsim$ mean.

For a measurable (we always assume it tacitly) set $E$ we denote by $\abs{E}$ its $n+1$ dimensional Lebesgue measure. For the integral average, we have the standard notation
\[f_E = \dashint_{E} f = \frac{1}{\abs{E}} \int _{E} f.\] 
For a (measurable) function $f$, we denote
\[f^{+} = f_{+} = f 1_{\{f > 0\}} \quad \textrm{and} \quad f^{-} = f_{-} = -f 1_{\{f < 0\}}. \]
The function $1_E$ equals $1$ in $E$ and zero elsewhere.

We continue by introducing the notation for parabolic rectangles. Let $p > 1$ be a number that is fixed throughout the paper. In applications, it would be the $p$ coming from the $p$-Laplace operator. For evolutionary problems, it has an important role in determining how the time variable and the space variables scale.

If $Q \subset \mathbb{R}^{n}$ is a cube with sides parallel to coordinate axes, we denote its side length by $\ell(Q)$. Take a parameter $\gamma \in (0,1)$. In previous papers, this has been called \textit{the lag}, but here the name \textit{shape} would be better. We specify a parabolic rectangle together with its upper and lower parts by its center $(x,t)$, side length $\ell(Q)$, and shape parameter $\gamma$. We usually drop some or all of the parameters from the notation, but if they are present, they should be understood as follows (see also Figure \ref{fig2:pararectangle}):
\begin{align*}
R((x,t),\ell(Q),\gamma) &= Q \times (t-\ell(Q)^{p},t+\ell(Q)^{p}) \\
R^{-}(\gamma) &= Q \times (t-\ell(Q)^{p},t-(1-\gamma))\ell(Q)^{p}) \\
R^{+}(\gamma) &= Q \times (t+(1-\gamma)\ell(Q)^{p},t+\ell(Q)^{p}).
\end{align*}
The number $\ell(R) :=  \ell(Q)$ is called the side length of the rectangle. Addition of a constant to a set in $\mathbb{R}^{n+1}$ is always understood as adding the constant to the time coordinate. 

\subsection{Classes of $\BMO$ type}
\label{subse:bmo}
We say that $u \in \PBMO^{-}$ if each parabolic rectangle $R$ has a constant $a_{R}$ such that 
\[\norm{u}_{\PBMO^{-}} := \sup_{R} \left( \dashint_{R^{-}(\frac{1}{2})} (u-a_R)^{+} + \dashint_{R^{+}(\frac{1}{2})}(a_R - u)^{+} \right) < \infty. \]
This is not a norm in the precise meaning of the word, but we call it a norm. $\PBMO^{-}$ is not a vector space, but we call it a space. Moreover, we say that an operator $T$ is bounded on $\PBMO^{-}$ if $\norm{Tu}_{\PBMO^{-}} \leq C \norm{u}_{\PBMO^{-}}$ even if the set up is not the one of normed linear spaces. 

The methods developed in the previous works show that given \textit{shape} $\gamma$ and \textit{lag} coefficient $L  > (1-\gamma)$, there are constants $\{b_R\}_{R}$ such that
\begin{equation*}
\sup_{R} \left( \dashint_{R^{-}(\gamma)} (u-b_R)^{+} + \dashint_{R^{-}(\gamma) + L \ell(R)^{p}}(b_R - u)^{+} \right) \eqsim_{n,p,\gamma,L} \norm{u}_{\PBMO^{-}}.
\end{equation*}
Suitable references for this are \cite{KS2016_1,Saari2014}, and the idea of the proof is also contained in the Lemma \ref{lemma:chain} proved in this paper. That lemma is intended to be an easy reference for the numerous applications of the chain argument that we need. 

It has been observed already earlier that $\PBMO^{-}$ can be realized as an intersection of two even rougher function classes. We mention \cite{KS2016_1} as a reference for the multidimensional case. In dimension one, this claim does not make so much sense because the ``rough'' function classes turn out to coincide and equal to $\PBMO^{-}$. See \cite{MRT1994}. However, when it comes to so called one-sided $L^{\infty}$ functions, similar things also happen in dimension one \cite{AC1998}. Let  
\begin{align*}
\norm{u}_{\BMO^{+}(\gamma,L)} := \sup_{R}  \dashint_{R^{-}(\gamma)} (u-u_{R^{-}(\gamma) + L \ell(R)^{p}})^{+} &< \infty \\
\norm{ -u }_{\BMO^{-}(\gamma,L)} := \sup_{R} \dashint_{R^{-}(\gamma) + L\ell(R)^{p}}(u_{R^{-}(\gamma)} - u)^{+} &< \infty.
\end{align*}
The first inequality is called $\BMO^{+}$ condition and the second one $ - \BMO^{-}$ condition. Note that $\BMO^{-}$ would just be $\BMO^{+}$ with $t$-axis of the coordinate space reversed. This notational convention on the direction of time also holds for $\PBMO^{\pm}$ and the maximal functions that we use. The \textit{one sided} function classes $\BMO^{\pm}$ are not known to be independent of $\gamma$ or $L$. However, they are useful because of the following:
\begin{equation}
\label{eq:pbmo_palasina}
\PBMO^{-} =  [\BMO^{+}(\gamma_{1},L_{1}) ] \cap [-\BMO^{-}(\gamma_{2},L_{2})]
\end{equation} 
for any choice of the shape and lag parameters. For details about this, see \cite{KS2016_1} and the discussion preceding the point to which we have advanced.

$\PBMO^{-}$ is closed under addition and multiplication by positive constants. Multiplication by negative constants reverses the direction of time, that is, it maps $\PBMO^{-}$ to $\PBMO^{+}$. For one-sided spaces $\BMO^{\pm}$ proving or disproving the previously mentioned property is an open problem.

We conclude this section by recalling the parabolic John-Nirenberg inequality. This appears in the literature, and it is proved in \cite{Aimar1988} whereas its formulation in different geometric configurations is studied in \cite{Saari2014}.
\begin{lemma}
\label{lemma:JN}
Let $u \in \PBMO^{-}$ and $\gamma \in (0,1)$. Take a parabolic rectangle $R$. Then there are constants $b_R$ and $c_1, c_2 \eqsim_{n,p,\gamma} 1$ such that
\begin{align*}
 \dashint_{R^{-}(\gamma)} \exp \left( \frac{c_1}{\norm{u}_{\PBMO^{-}}} (u-b_R)^{+} \right)  \leq c_2 \\
 \dashint_{R^{+}(\gamma)} \exp \left( \frac{c_1}{\norm{u}_{\PBMO^{-}}} (b_R-u)^{+} \right)  \leq c_2.
\end{align*}
\end{lemma}

\subsection{Parabolic maximal function}
\begin{figure}
\includegraphics[scale=0.7]{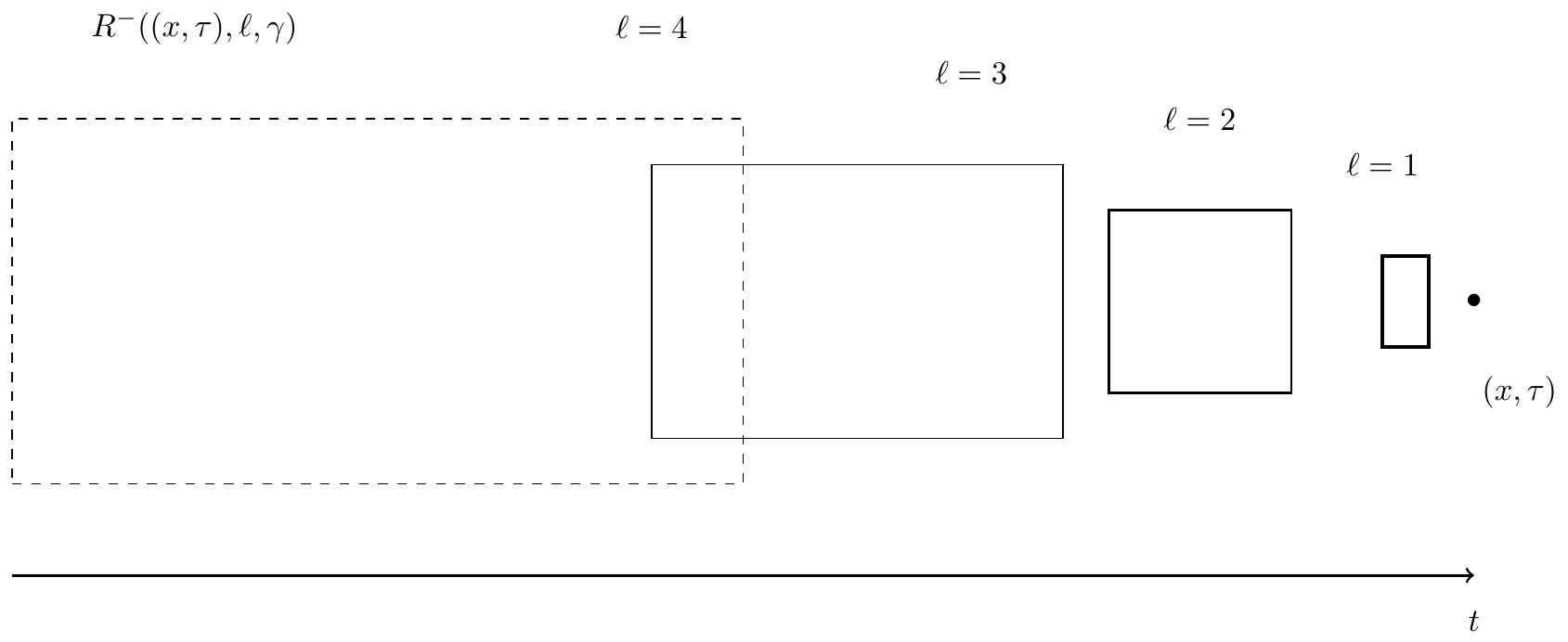}
\caption{\label{fig3:maximal_function}Sketch of the region in which the maximal function sees the positive part of the function.} 
\end{figure}

The first candidate to be a parabolic backward-in-time maximal operator was introduced in \cite{KS2014}, and it reads as
\[M^{\gamma-}f(x) = \sup_{R(x)} \dashint_{R^{-}(\gamma)} \abs{f}.\]
The supremum is over parabolic rectangles centred at $x$ and we average the absolute value over the left part. See also Figure \ref{fig3:maximal_function}. This definition is problematic when dealing $\PBMO^{-}$ functions that are not necessarily positive. Too many of them are mapped to infinity. For the positive functions, however, it coincides with the following maximal function.

We let 
\begin{align*}
M_*^{\gamma-} f (x)&:= \sup_{\ell > 0 } \left((f^{+})_{R^{-}(x,\ell,\gamma)} +  (f^{-})_{R^{+}(x,\ell, \gamma)}\right) \\
&:= \sup_{R(x)} \left( \dashint_{R^{-}(\gamma)} f1_{\{f > 0\}} - \dashint_{R^{+}(\gamma)} f 1_{\{f < 0\}} \right) 
\end{align*}
where the supremum is over all the parabolic rectangles that are centred at $x$. This maximal function measures \textit{positivity} in the past, \textit{negativity} in the future. In case the time variable is trivial, that is, we have functions on $\mathbb{R}^{n}$ as in the study of elliptic PDE or as in the classical Calder\'on-Zygmund theory, it coincides with the usual Hardy-Littlewood maximal function. 

There are, however, delicate issues involved in the behaviour of positive and negative parts of the functions, multiplications by negative numbers, and taking absolute values in the parabolic theory. Recall for instance the problem of $- \BMO^{-}$, and see Remark \ref{remark:katkonta}.  
In connection with $\PBMO^{-}$, we have that $\PBMO^{-} = - \PBMO^{+}$. The maximal function above is designed to have the same property while maintaining the consistency with the classical Hardy-Littlewood maximal function. Namely, $M_*^{\gamma-}(- f) = M_*^{\gamma +} f$. Note also that this maximal function has its natural \textit{sharp} version that can be used to define $\PBMO^{-}$. 

\section{Preliminary results}
\label{sec:preli}
This section contains the proofs of some lemmas that we need. In many occasions, we have to decompose functions in $\PBMO^{-}$ into positive and negative parts, and we next prove that this procedure respects the $\PBMO^{-}$ property. The subsequent subsection contains a chain lemma, and the remaining two contain a construction of a dyadic grid suitable for our purposes and some remarks on the local integrability of maximal functions of functions in $\PBMO^{-}$. 




\subsection{Truncations}
We start with noting that $\PBMO^{\pm}$ classes are stable with respect to truncation. Even if this property looks very elementary, its proof still differs a bit from the classical analogue. Namely, we really need the equivalence with respect to lag and shape in order to get the conclusion.

\begin{proposition}
\label{prop:katkot}
Let $u \in \PBMO^{-}$. Then $u^{+}, -u^{-} \in \PBMO^{-}$.
\end{proposition}
\begin{proof}
Take $u \in \PBMO^{-}$. We start with $u^{+}$ and the $\BMO^{+}$ bound. Take any $\gamma \in (0,1)$. By elementary considerations
\begin{align*}
\dashint_{R^{-}(\gamma)}(u^{+} - (u^{+})_{R^{+}(\gamma)})^{+} & \leq \dashint_{R^{-}(\gamma)}(u^{+} - (u_{R^{+}(\gamma)})^{+})^{+} \\
 & \leq \dashint_{R^{-}(\gamma)}(u- u_{R^{+}(\gamma)})^{+} \lesssim \norm{u}_{\PBMO^{-}}. 
\end{align*}
For the $\BMO^{-}$ bound, take a parabolic rectangle $R$ with side length $\ell(R)$. Let
\begin{align*}
Q = R^{-}(\gamma), &\quad Q^{+} = R^{+}(\gamma) \quad \textrm{and}\\
 & Q^{++}= R^{+}(\gamma) + (0, \ldots, 0, (1+\gamma) \ell(R)^{p}).
\end{align*}
Then it is easy to see that for every point in $Q^{++}$ we have that
\begin{align*}
( (u^{+})_{Q} - u^{+})^{+} &\leq ([(u-u_{Q^{+}})^{+}]_{Q} + (u_{Q^{+}})^{+} - u^{+} )^{+} \\
 &\leq  \dashint_{Q}(u-u_{Q^{+}})^{+} + (u_{Q^{+}} - u)^{+} \\
 &\lesssim \norm{u}_{\PBMO^{-}} + (u_{Q^{+}} - u)^{+},
\end{align*}
and averaging over the domain $Q^{++}$ completes the proof for $u^{+}$. See Remark \ref{remark:ketjun_yleistys} for the fact that the increased gap between $Q$ and $Q^{++}$ does not matter. 

For the negative part, we note that $-u \in \PBMO^{+}$ and by the previous step also $(-u)^{+} \in \PBMO^{+}$. Then $-u^{-} = - (-u^{+}) \in - \PBMO^{+} = \PBMO^{-}$.
\end{proof}

\begin{remark}
\label{remark:katkonta}
The previous proposition does not say anything about absolute values of $\PBMO^{-}$ functions. Indeed, $\PBMO^{-}$ is not closed under subtraction. In general, absolute values may fail to belong to $\PBMO^{-}$. An easy one dimensional example is $e^{t} - e^{-t}$. As an increasing function it trivially satisfies $\PBMO^{-}$ condition. However, taking an absolute value converts it into a rapidly decreasing function as we go towards $- \infty$.
\end{remark}


\subsection{Chain argument}
\label{subse:chain}
The following lemma is what we refer to as the chain argument. It is very important. The formulation and the proof simplify the arguments that were used in \cite{Saari2014} a lot, and we hope that writing this argument down in a readable way is of independent interest. Moreover, as a corollary of this lemma, we can also easily deduce the independence of shape and lag of the definition of $\PBMO^{-}$ which we often refer to.

\begin{lemma}
\label{lemma:chain}
Let $u \in \PBMO^{-}$. Consider a rectangle $R^{-}(\theta)$ with spatial side length $\ell$. Let $v \in \mathbb{R}^{n}$ and $\tau > (1-\theta)$. Define 
\[P^{+}(\theta) = R^{-}(\theta) + (v,\tau \ell ^{p}). \]
Then
\[ \dashint_{P^{+}(\theta)}  \dashint_{R^{-}(\theta) } (u(x)   - u(y) )^{+} \dx \dy \lesssim_{n,p,\theta,\abs{v}/\ell, \tau} \norm{u}_{\PBMO^{-}}       \]
\end{lemma}
\begin{proof}
This is a very simple instance of the chain argument, but for the reader's convenience we give a proof. Assume first that $\tau$ is very large. Denote the base of $R^{-}(\theta)$ by $Q$. It has side length $\ell$. Let 
\[Q_i = Q + \frac{v}{\abs{v}} r_i, \quad r_i > 0 . \]
We choose such a sequence $\{r_i \}_{i=1}^{k}$ that
\[\frac{\abs{ Q_i \cap Q_{i+1} }}{\abs{Q_0}} = \delta \]
where $\delta \in (0,1)$ is a number such that $Q_k = Q_0 + v$. Note that the numbers $\delta$ and $k$ depend only on $\abs{v}$.

\begin{center}
\begin{figure}
\includegraphics[scale=0.8]{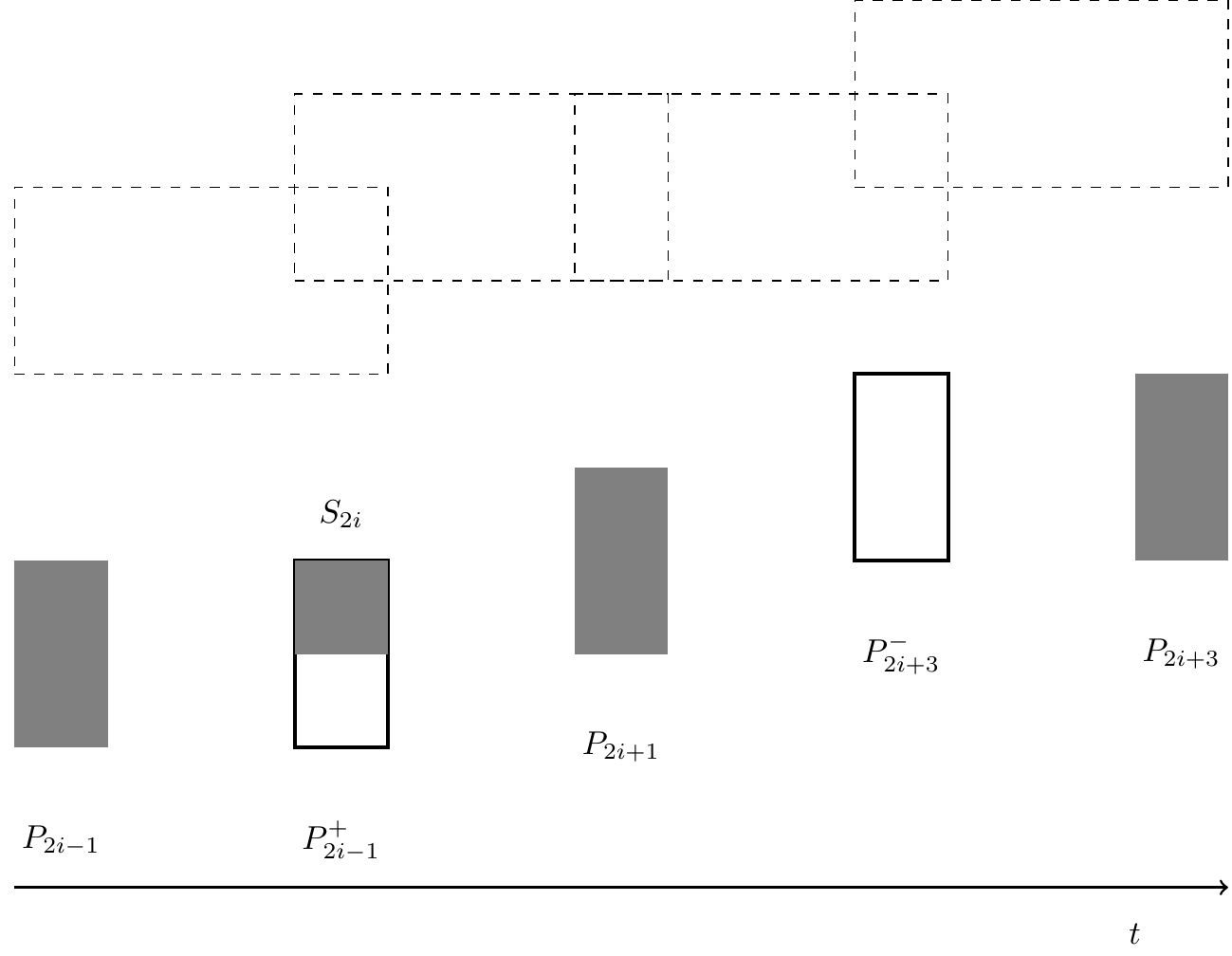}
\caption{\label{fig1:chain}The sets $P_{2i-1}$, their translates forwards and backwards in time, and the real parabolic rectangles motivating these boxes.} 
\end{figure}
\end{center}

Let $P_1 = R^{-}(\theta)$. It has the least time coordinate $t_1$. Define
\begin{align*}
P_{2i - 1} &= Q_{i} \times (t_1 + (i-1) (1 + \theta)\ell^{p},t_1 + [(i-1) (1 + \theta) + (1-\theta) ]\ell^{p} ) \\
P_{2i-1}^{+} &= P_{2i - 1} + (1+ \theta) \ell^{p}  \\
P_{2i-1}^{-} &= P_{2i - 1} - (1+ \theta) \ell^{p}\\
S_{2i} &= P_{2i-1}^{+} \cap P_{2i+1}^{-}  .
\end{align*}
The sets $S_{2i}$ are such that $\frac{\abs{ S_{2i} }}{\abs{P_{2i-1}}} =\delta$. Moreover, the pairs of sets $(P_{2i - 1},P_{2i - 1}^{+})$ can be used in testing the $\PBMO^{-}$ conditions. See also the Figure \ref{fig1:chain}.

Since the chain of cubes travels the amount $\abs{v}$ in space, the number of cubes needed is roughly 
\[k \eqsim \abs{v} / \ell.\] 
This means that the sequence of space time rectangles $\{P_{2i-1}\}_i$ advances, up to a constant factor, the amount $\abs{v} \ell^{p} / \ell = \abs{v} \ell^{p-1}$ in time.

By our assumption on $\tau$ being large, we can find $\delta$ such that 
\[P_{2(k +l)+1 } = P^{+}(\theta)\] 
where $l \lesssim \tau$ is an integer taking care of the possibility that $\tau$ is too large. The exact requirement on $\tau$ is 
\begin{equation}
\label{eq:tau_raja}
\tau \ell^{p} \gtrsim \abs{v} \ell^{p-1}.
\end{equation}

Note first that
\begin{align*}
&\dashint_{P^{+}(\theta)} \dashint_{R^{-}(\theta) } (u(x)   - u(y) )^{+} \dx \dy \\
&\leq \dashint_{R^{-}(\theta)} (u  - u_{R^{+}(\theta)}  )^{+} 
+ (u_{R^{+}(\theta)} - u_{P^{-}(\theta)})^{+}
+ \dashint_{P^{+}(\theta)} (u_{P^{-}(\theta)}  - u  )^{+} \\
&\lesssim (u_{R^{+}(\theta)} - u_{P^{-}(\theta)})^{+} + 2\norm{u}_{\PBMO^{-}}.   
\end{align*}
By this reduction, there is no loss in generality in the following estimation:
\begin{align*}
(u_{R^{-}(\theta)}  &- u_{P^{+}(\theta)}  )^{+} \leq \sum_{i=1}^{k+l} (u_{P_{2i-1}} - u_{P_{2i+1}})^{+} = \sum_{i=1}^{k+l} \dashint_{S_{2i}}  (u_{P_{2i-1}} - u_{P_{2i+1}})^{+} \\
& \leq  \sum_{i=1}^{k+l}  \left( \dashint_{S_{2i}}  (u_{P_{2i-1}}- u)^{+}+ \dashint_{S_{2i}}( u - u_{P_{2i+1}})^{+} \right) \\
& \lesssim_{\delta} \sum_{i=1}^{k+l}  \left( \dashint_{P_{2i-1}^{+}}  (u_{P_{2i-1}}- u)^{+}+ \dashint_{P_{2i+1}^{-}}( u - u_{P_{2i+1}})^{+} \right)\\
& \lesssim_{k,l} \norm{u}_{\PBMO^{-}}.
\end{align*} 
In this estimate, the quantities $l,k,\delta$ depend only on $\abs{v}/ \ell$, $\tau$, $n$, $p$ and $\theta$. 

To get rid of the restriction of $\tau $ being a large number, we do the following. We find an auxiliary pair of rectangles $P_\epsilon^{+}(\theta)$ and $R_{\epsilon}^{-}(\theta)$ that have side lengths $\epsilon \ell$ and that form upper and lower halves of a parabolic rectangle $R_{\epsilon}$ that is centred in the middle of the line connecting the centres of $P^{+}(\theta)$ and $R^{-}(\theta)$. 

We cut the original blocks into pieces with the same side length $\epsilon \ell$. The pieces of $R^{-}(\theta)$ are called $A_i^{-}$ and the pieces of $P^{+}(\theta)$ are called $B_i^{+}$. Then
\[(u_{A_i^{-}} - u_{B_j^{+}})^{+} \leq (u_{A_i^{-}} - u_{R_{\epsilon}^{-}(\theta)})^{+}+(u_{P_\epsilon^{+}(\theta)}- u_{B_j^{+}})^{+} .\]
These two terms are back in the original situation with data $v_i,v_j \in \mathbb{R}^{n}$, $\abs{v_i} \eqsim \abs{v_j} \eqsim \abs{v}$, and  $\tau_i \eqsim \tau_j \eqsim \frac{\tau}{\epsilon^{p}}$. Then
\[\tau_i (\epsilon \ell)^{p} \gtrsim \tau \ell^{p} ,\]
so the lower bound asked by \eqref{eq:tau_raja} is satisfied if
\[\tau \ell^{p} \gtrsim \abs{v_i} (\epsilon \ell) ^{p-1} \eqsim \frac{\abs{v}}{\ell} \epsilon^{p-1} \ell^{p},\]
whence we deduce that $\epsilon \lesssim ( \ell \tau /\abs{v})^{1/ (p-1)}$ can be made sufficiently small. This bound only depends on legal quantities, and we may consider $\epsilon$ a constant in the future.

Finally, we may compute
\begin{align*}
\dashint_{P^{+}(\theta)}  \dashint_{R^{-}(\theta)} (u(x)  - u(y) )^{+} \dx \dy &\lesssim  \sum_{i,j} \dashint_{{A_i^{-}}}\dashint_{{B_i^{-}}} (u(x) - u(y))^{+}
\end{align*}  
which reduces the situation to the case already handled.
\end{proof}

\begin{remark}
\label{remark:ketjun_yleistys}
Given two, possibly rather weird, bounded Borel sets $A$ and $B$ from $\mathbb{R}^{n+1}$ such that 
\[ \Delta =   \inf \{t: (x,t) \in B \} - \sup\{t: (x,t) \in A \} > 0 ,\]
we can conclude by the previous lemma that 
\[(u_A - u_B)^{+} \leq  \dashint_{A} \dashint_{B} (u(x) - u(y))^{+} \dx \dy \lesssim \norm{u}_{\PBMO^{-}} \]
where the dependency is on the parameters of the Lemma \ref{lemma:chain} associated to a pair of parabolic rectangles $R_1$ and $R_{2}$ that contain the sets $A$ and $B$ and satisfy the assumptions of the Lemma. More precisely, there is dependency in terms of a function separately increasing in $\abs{R_1}/ \min\{\abs{A}, \abs{B}\}$, $\max \{d(A),d(B)\} / \Delta $, and $d(A,B)/ \Delta$. 
\end{remark}

\begin{remark}
\label{remark:john-nirenberg}
The idea of the previous remark also generalizes to the John-Nirenberg inequality of $\PBMO^{-}$, Lemma \ref{lemma:JN}. This is quite a direct consequence of the convexity of the exponential. Hence we may give an improved formulation of the John-Nirenberg inequality with $c_1$ and $c_{2}$ depending on the same constants as above:
\[\dashint_{A} \dashint_{B} \exp\left( \frac{c_1}{\norm{u}_{\PBMO^{-}}} (u(x) - u(y))^{+} \right) \dx \dy \lesssim c_2.\]
\end{remark}


\subsection{Dyadic grid}
In course of the proof of the main theorem, we decompose a function $u$ into a bounded part and into additional pieces that are roughly of the form
\[b_Q = 1_{Q}(u - u_{Q^{+}} ) .\]
This form is particularly convenient when working with $u \in \PBMO^{-}$. We want the pieces to have supports with controlled overlap, and obtaining that property in the geometry where the time variable scales as space to power $p$ requires some work.

For an effective use of the Calder\'on-Zygmund stopping time argument, it is beneficial to dispose of a dyadic grid. The existence of such a grid is nontrivial in the parabolic geometry if we want to maintain an intuition of how our dyadic parabolic rectangles look like. This intuition is lost when using the black boxes that would be available from the analysis in metric spaces.

The main problem in the dyadic grid is that if the number $p$ governing the geometry is an irrational number, then we cannot simply subdivide a box to $2^{n}$ cubes in space and $2^{p}$ boxes in the space since the latter number is not an integer. This problem is easy to circumvent if $p$ is rational, but in case of irrational $p$ one has to resort to an approximative construction.

\begin{lemma}
\label{lemma:dyadic}
Let $Q_0$ be a half of a parabolic rectangle with side length $\ell(Q_0)$. Then there exists boxes with the properties of the dyadic tree that are almost parabolic sub-rectangles of $Q_0$:
\begin{enumerate}
\item[(i)] $\mathscr{D} = \cup_{i=0}^{\infty} \mathscr{D}_i$ and $\mathscr{D}_i = \{Q_j^{(i)} \}_{j}$. In addition $\cup_{j} Q_j^{(i)} = Q_0 \in \s{D}_{0}$ for all $i$. The boxes $Q_j^{(i)}$ with common $i$ are translates of each other. 
\item[(ii)] If $P, Q \in \s{D}$, then $P \cap Q \in \{\varnothing, P, Q\}$. For each $Q \in \s{D}_i$, there is a unique $\widehat{Q} \in \s{D}_{i-1}$ with $\widehat{Q} \supset Q$. In addition, $\abs{Q} \eqsim \abs{\widehat{Q}}$.
\item[(iii)] Every  $Q \in \s{D}_{i}$ almost has the dimensions $(\ell, \ell^{p})$. Namely, there is a parabolic box $\tilde{Q}$ obtained from $Q_0$ by means of parabolic dilation $(x,t) \mapsto (\delta x, \delta^{p}t)$ and translation so that $\ell(\tilde{Q}) = 2^{-i} \ell(Q_0)$, $\tilde{Q} \supset Q$, and $\abs{\tilde{Q}} \eqsim \abs{Q}$
\end{enumerate}
\end{lemma}
\begin{proof}
For simplicity, and without loss of generality, we may assume that $Q_0 = [0,1]^{n+1}$.
For every $p > 1$, it is easy to see that there exists a non-decreasing sequence of integers $k_i$ such that for $i \in \mathbb{Z}_{+}$
\[\left \lvert p - \frac{k_i}{i} \right \rvert \leq \frac{1}{i} .\]
We denote $q_i = k_i/i$.
  
At step one, we divide the spatial side length of $Q_0$ by $2$ and the temporal one by $2^{k_1}$. These dimensions give us space time boxes in the geometry where time scales as space to power $q_1$. We use them to partition $Q_0$. This gives the collection $\mathscr{D}_1$ that serves as a generation in the dyadic grid we are constructing. 

At the second step, we keep repeating the process so that boxes in $\mathscr{D}_{i-1}$ have their spatial side length halved so that we reach the desired side length $2^{-i}$. The previous side length was $2^{1-i}$. We choose the temporal side length to be $2^{-k_i}$. Note that this is obtained by repeatedly multiplying $2^{-k_{i-1}}$ by $2^{-1}$ (possibly zero times) since the sequence $k_i$ is increasing. With these dimensions, we get space time boxes in the geometry where time scales as the space to power $q_i$. We use these blocks to partition boxes in $\s{D}_{i-1}$. Their collection is called $\s{D}_{i}$. 

Finally note that if $Q \in \s{D}_i$, and $\tilde{Q}$ is a block with the same spatial side length that respects the geometry where time scales as space to power $p$, then 
\begin{align*}
\frac{\abs{Q}}{\abs{\tilde{Q}}} = \frac{2^{-k_i}}{2^{-p i}} = 2^{p i - k_i} \in (2^{-1},2).
\end{align*}
Hence we have good control over the distortion of the geometry at all scales. Consequently, we may always cover the dyadic box with a proper parabolic rectangle so that the error in volume is controlled.
\end{proof}

\subsection{Integrability of maximal functions}
The functions in $\PBMO^{-}$ are locally integrable. In spatially unbounded domains the maximal functions may easily fail this property so some care must be taken. One suitable criterion for local integrability of maximal function of a $\PBMO^{-}$ function is given by its finiteness.

\begin{lemma}
\label{lemma:finite-Lloc}
\begin{enumerate}
\item[(i)] Let $u \in \PBMO^{-}(\Omega)$ where $\Omega \subset \mathbb{R}^{n+1}$ is a space time cylinder. Then $M_*^{\gamma-}u \in L_{loc}^{1}(\Omega)$. 
\item[(ii)] Let $u \in \PBMO^{-}(\mathbb{R}^{n+1})$. If there are $x_1,x_2 \in \mathbb{R}^{n}$ and numbers $t_2 > t_1$ such that
\begin{align*}
M^{\gamma-} u^{+}(x_2,t_2) &< \infty \quad \textrm{and} \\
M^{\gamma+} u^{-}(x_1,t_1) &< \infty,
\end{align*}
then 
\[M_*^{\gamma-} u \in L_{loc}^{1}(\{(y,\tau) \in \mathbb{R}^{n+1}: t_1< \tau < t_2 \} ).\]
\end{enumerate}
\end{lemma}
\begin{proof}
(i) For the case of the space time cylinder, choose any rectangle $R \subset \Omega$. We may apply the John-Nirenberg from \cite{Saari2014} to get exponential integrability for $u^{\pm}$ in the union of the rectangles that the maximal function sees from there. This is enough to show that the maximal function has to be locally integrable.

(ii) In the case where $u$ is defined on $\mathbb{R}^{n+1}$, some more analysis is needed. Assume first that $u \geq 0$. We may restrict our attention to a maximal operator that only sees large rectangles since the case of small rectangles is the same as the case of a space time cylinder. If we write $U = \max\{U_1, U_2 \}$ where $U$ is the maximal function, $U_1$ is the small-rectangle-supremum, and $U_2$ is the supremum over large rectangles, we see that for all parabolic rectangles $R$ located in the half-space 
$ \{(y,\tau) \in \mathbb{R}^{n+1} : \tau < t_2\} $
we have that
\[\dashint_{R^{-}(\gamma)} U_2 \lesssim_{R,n,p,\norm{u}_{\PBMO^{-}}} U(x,t) < \infty\]
by the argument in Section \ref{sect:large-rectangle}. Hence we get that $M_*^{\gamma-} u$ is locally integrable in the half-space. In case we have a general $u = u^{+} - u^{-}$, we may do the previous estimation for the positive and negative parts separately.
\end{proof}


\section{Parabolic BMO and maximal functions}
\label{sec:para}
Next we study the boundedness of the parabolic maximal operator on parabolic $\BMO$ space. This is inspired by a result in Bennett, DeVore, and Sharpley \cite{BdVS1981}. Indeed, under the positivity assumption, the following lemma generalizes their result. The case of more general functions is discussed after the proof of this lemma.

We do not distinguish the cases where the domain of definition of $\PBMO^{-}$ functions is a space time cylinder or a full space. The only difference between these cases is in integrability issues of the maximal function. In what follows, we assume the local integrability of $M_*^{\gamma-} u$, and criteria for this condition in the two different cases can be found in Lemma \ref{lemma:finite-Lloc}.

\begin{lemma}
\label{lemma:boundedness}
Let $u \in \PBMO^{-}$ be such that $u \geq 0$ almost everywhere. Then $M^{\gamma-}u = M_*^{\gamma-}u \in \PBMO^{-}$ and 
\[\norm{M_*^{\gamma-}u}_{\PBMO^{-}} \lesssim_{\gamma} \norm{u}_{\PBMO^{-}} \]
provided that $M_*^{\gamma-}u$ is locally integrable. 
\end{lemma}
\begin{proof}
Let $u \geq 0$ be in $\PBMO^{-}$. We start by proving that $M_*^{\gamma-} u \in \BMO^{+}$. Let $R_0$ be an arbitrary parabolic rectangle. We will use the following notation:
\begin{align*}
U(x) &= M_*^{\gamma-}u(x) \\
U_1(x) &= \sup \{u_{R^{-}(\gamma)}: \textrm{$x$ is the center of $R$},  \ell(R) \leq 100^{-1} \ell(R_0)   \} \\
U_2(x) &= \sup \{u_{R^{-}(\gamma)}: \textrm{$x$ is the center of $R$}, \ell(R) \leq 100^{-1} \ell(R_0)  \}. 
\end{align*}
Note that $U(x) = \max \{U_1(x),U_2(x)\}$. We may choose the configuration in which we attempt to bound the $\PBMO^{-}$ norm type quantities. This is due to Lemma \ref{lemma:chain}. We let $R_0^{-} = R_0^{-}(0)$ be a standard half of a parabolic rectangle whereas $R_0^{+} = R_0^{-} + 100 \ell^{p}(R_0^{-})$.

We will obtain estimates for 
\[\dashint_{R_0^{-}}(U_1(x) - U_{R_0^{+}})^{+} \quad \textrm{and} \quad \dashint_{R_0^{-}}(U_2(x) - U_{R_0^{+}})^{+}.\] 
The proof works for a much more general expression, but in order to emphasize the structure of the proof, we try to avoid excessive technicalities. We start the proof with the case of $U_2$.

\subsection{Large rectangles.}
\label{sect:large-rectangle}
Take a point $x \in R_0^{-}$ and any rectangle $R $ that is admissible in the definition of $U_2$, that is, $x$ is the center of $R$ and $\ell(R) \geq 100^{-1} \ell(R_0)$. We denote by $A$ the number satisfying $\ell(R) = A \ell(R_0) $  

Take a point $z$ from $R^{+}$ and a rectangle $P$ centered at $z$ with
\[u_{P^{-}(\gamma)} \leq M_*^{\gamma-}u(z)  \]
such that it has large enough side length $B \ell(R_0)$. The constant $B$ will be fixed later. It has to have the correct relations to the constant $A$.

Next we look at how close the rectangles $R^{-}(\gamma)$ and $P^{-}(\gamma)$ are to each other. Let 
\begin{align*}
\tau^{-} &= \sup \{t: (x,t) \in R^{-}(\gamma)  \} \\
\tau^{+} &= \inf \{t: (x,t) \in P^{-}(\gamma)  \}.
\end{align*}   
Note that since
\[\tau^{+} - \tau^{-} \geq \ell(R_0)^{p}(A^{p}\gamma - B^{p}), \]
we may choose the constant $B$ to satisfy
\[B^{p} = \frac{\gamma}{100} A^{p} \] 
so that $P$ and $R$ are of comparable size and well separated. Independently of $R$, the separation of the rectangles in both time and space variables is governed by the parameters $n,p,\gamma$. 
Hence me may apply the Lemma \ref{lemma:chain} in form of Remark \ref{remark:ketjun_yleistys}:
\begin{align*}
( u_{R^{-}(\gamma)}- U(z) )^{+} & \leq    ( u_{R^{-}(\gamma)} - u_{P^{-}(\gamma)} )^{+} \nonumber \\
&\lesssim \norm{u}_{\PBMO^{-}} .
\end{align*}
Taking a supremum over $R$, we get a pointwise bound for $U_2- U(z)$ for every $z \in R_0^{+}$. This also proves the claim about the mean value:
\begin{equation}
\label{eq:large_double_control}
 \dashint_{R_0^{+}} \dashint_{R_0^{-}}  (U_2(x)  - U(y) )^{+} \dx \dy \lesssim \norm{u}_{\PBMO^{-}}.
\end{equation}
This estimate controls the large rectangle part trivially. 

\begin{remark}
Note that in this part it was only essential that the positive term was average of a large rectangle, and the size of the rectangle in the negative term could be chosen. Hence we got
\[           \sup_{r > 100^{-1} \ell(R_0)} \inf_{\rho > 0}   ( u_{R^{-}(x,r,\gamma)} - u_{R^{-}(y,\rho,\gamma )})^{+} \lesssim \norm{u}_{\PBMO^{-}}. \]
\end{remark}

\subsection{Small rectangles I: Calder\'on-Zygmund.}
\label{sec:CZ}
In this part, we use the dyadic grid constructed in Lemma \ref{lemma:dyadic}. The cases with $p \in \mathbb{Q}$ and $p \notin \mathbb{Q}$ are not very different since the fact that the dyadic grid only approximates the real one only comes into the picture in few occasions. In order to keep the notation more simple, we work on the case with rational $p$, and comment the corrections that should be done with approximate dyadic grid. This saves us from some additional indices. Before attacking our target rectangle $R_0$, we make some preparations.

Take a parabolic rectangle with lower part $Q_0^{-}$. Let $\Omega_1 = \{x \in Q_0^{-} : u_{Q_0^{+}} < U_1(x) \}$. We run a Calder\'on-Zygmund stopping time argument on $Q_0^{-}$ with a stopping rule $u_{Q_i^{+}} > \lambda$ where $\lambda > u_{Q_0^{+}}$. Denote by $\widehat{Q}$ the dyadic parent of $Q$. We decompose $u = b+g$ where the components are
\begin{align*}
b_i & = 1_{Q_i} (u-u_{\widehat{Q_i}^{+}}), \\
b &= \sum_i b_i, \\
g &= \sum_{i} 1_{Q_i} u_{\widehat{Q_i}^{+}}  + u1_{Q_0^{-} \setminus \cup_{i} Q_i}.
\end{align*}
In case we use the approximate dyadic grid, the stopping rule applies to $\left(\tilde{\widehat{Q}} \right)^{+}$ rectangles while the indicator functions in the decomposition of $u$ are left as they are.  Recall the details of the dyadic grid from Lemma \ref{lemma:dyadic}.

By construction, $\norm{g}_{L^{\infty}} \leq \lambda$. By the parabolic John-Nirenberg inequality (Lemma \ref{lemma:JN} and Remark \ref{remark:john-nirenberg}), we have that 
\[\dashint_{\widehat{Q_i}} e^{2 \epsilon b_i^{+}} = \dashint_{\widehat{Q_i}} e^{ 2 \epsilon (u-u_{\widehat{Q_i}^{+}}) } \leq 2\]
if $\epsilon \lesssim \norm{u}_{\PBMO^{-}}^{-1}$ is small enough. By elementary properties of the maximal functions
\begin{align*}
\dashint_{Q_0^{-}} e^{\epsilon M_*^{\gamma-} b^{+}} & \leq \left( \dashint_{Q_0^{-}}  (M_*^{\gamma-} e^{\epsilon b^{+}} )^{2} \right)^{1/2} 
\lesssim  \left( \dashint_{Q_0^{-}}  e^{ 2  \epsilon b^{+} } \right)^{1/2}.
\end{align*}
By the fact that $Q_i$, the supports of $b_i,$ are disjoint, we conclude that 
\begin{align*}
 \dashint_{Q_0^{-}}  e^{ 2  \epsilon b^{+} } & \leq 1 +  \sum_{i}  \frac{\abs{Q_i}}{\abs{Q_0^{-}}} \dashint_{Q_i} e^{2 \epsilon b_i} \lesssim  1+ \sum_{i}  \frac{\abs{Q_i}}{\abs{Q_0^{-}}} \dashint_{\widehat{Q_i}} e^{2 \epsilon b_i} \lesssim 1.
\end{align*}

Moreover,
\[ M_*^{\gamma- } (1_{Q_0^{-}} u)= M_*^{\gamma- } (1_{Q_0^{-}} (b+g)^{+})  \leq M_*^{\gamma- } (b^{+}) + M_*^{\gamma- } (g^{+}) \]
so using the previous computation, we obtain
\begin{align*}
\int_{ \Omega_1} U_1(x) &\lesssim  \abs{Q_0^{-}} \norm{u}_{\PBMO^{-}} \log \exp \dashint_{Q_0}   M_*^{\gamma-}( \epsilon b^{+})  +   \abs{\Omega_1} \norm{M_*^{\gamma-} (g^{+})}_{L^{\infty}} \\
&\leq C_0 \norm{u}_{\PBMO^{-}} \abs{ Q_0^{-}} +   \abs{\Omega_1} \norm{g^{+}}_{L^{\infty}} \\
&\leq \norm{u}_{\PBMO^{-}} \abs{ Q_0^{-} } + \lambda \abs{\Omega_1}.
\end{align*}
Subtracting $\abs{\Omega_1} \lambda $ and letting $\lambda \to u_{Q_0^{+}(\gamma)}$, we get
\[\dashint_{Q_0^{-}}(U_1(x) - u_{Q_0^{+}})^{+} \lesssim \norm{u}_{\PBMO^{-}}. \]

\subsection{Small rectangles II: Chains.}
Next consider the function $U_1$ in $R_0^{-}$. For computing the averages, the maximal function only has the region 
\[E = \bigcup_{\substack{z \in R_{0}^{-}(\gamma) \\ 0<r< 100^{-1}\ell(R_0)}}  R(z,r,\gamma)  \]
at its disposal. This region is small compared to the gap between $R_0^{-}$ and $R_0^{+}$, so we may cover it with one single rectangle $Q_0^{-}$ such that $Q_0^{+}$ is still before $R_0^{+}$.

Now we can localize the maximal function to the region $Q_0^{-}$.
For all $x \in R_0^{-}$, we have that
\[U_1 \leq M_*^{\gamma-}(1_{E} u) \leq M_*^{\gamma-}(1_{Q_0^{-}} u) .\]
This together with the previous Calder\'on-Zygmund consideration gives the estimate
\begin{align}
\label{eq:small_rectangle_final}
\int_{R_0^{-}} ( U_1 - U_{R_0^{+}} )^{+} & \leq \int_{Q_0^{-}} (M_*^{\gamma- }(1_{Q_0^{-}}u )- u_{Q_0^{+}})^{+}   + (u_{Q_0^{+}} - u_{R_0^{+} } )^{+} \nonumber \\
& \lesssim \abs{Q_0^{-}} \norm{u}_{\PBMO^{-}} + (u_{Q_0^{+}} - u_{R_0^{+} } )^{+} .
\end{align}
In the first term, recall that $\abs{Q_0^{-}} \eqsim \abs{R_0^{-}}$. For the second term, we may use the Remark \ref{remark:ketjun_yleistys} generalizing Lemma \ref{lemma:chain}. This completes the proof of the case of small rectangles.

\subsection{The case of all sizes.}
In general, we may estimate the full maximal function 
\[U = \max \{ U_1, U_2 \}\] 
by the sum of its two parts:
\begin{align*}
\int_{R_{0}^{-}}(U -U_{R_0^{+}})^{+} &\leq \int_{R_{0}^{-}\cap \{U_1 \geq U_2 \}}(U_1 -U_{R_0^{+}})^{+} \\
& \quad  + \int_{R_{0}^{-} \cap \{U_2 > U_1\}}(U_2 -U_{R_0^{+}})^{+} .
\end{align*} 
The first term was bounded in the previous subsections, and the second one was the large rectangle case.

Up to now we have proved that if $u \geq 0$ and $u \in \PBMO^{-}$, then $M_*^{\gamma-} u \in \BMO^{+}$. In order to get the claimed $M_*^{\gamma-}u \in \PBMO^{-}$, we still have to show that 
$U = M_*^{\gamma- }u \in -\BMO^{-}$, that is,
\begin{equation*}
\sup_{R} \dashint_{R^{+}} (U_{R^{-}} - U )^{+} \lesssim \norm{u}_{\PBMO^{-}}.
\end{equation*} 
Compare to the identity \eqref{eq:pbmo_palasina}.

However, this can be reduced to the case just handled. Take three blocks similar to $R_0^{\pm}$ in the previous considerations:
\begin{align*}
&Q = R_0^{-}  , \quad Q^{+} = R_0^{+}  \quad \textrm{and} \\
&Q^{++}= R_0^{+} + (0, \ldots, 0,  100 \ell(R_0)^{p}). 
\end{align*}
Now 
\begin{align*}
\dashint_{Q^{++}} &(U_{Q} - U )^{+} \leq \dashint_{Q^{++}} \dashint_{Q} ( U(x) - U(y) )^{+} \dx \dy \\
& =  \dashint_{Q^{++}} \frac{1}{\abs{Q}} \int_{Q\cap \{U_1 >  U_2 \}} ( U_1(x) - U(y) )^{+} \dx \dy  \\
& \qquad + \dashint_{Q^{++}} \frac{1}{\abs{Q}} \int_{Q\cap \{U_1 \leq  U_2 \}} ( U_2(x) - U(y) )^{+} \dx \dy \\
& = I + II.
\end{align*}
The second term $II$ is clear by the same argument that lead to the large rectangle case inequality \eqref{eq:large_double_control}. For the first term, we can estimate
\begin{align*}
I &\leq \dashint_{Q^{++}} \frac{1}{\abs{Q}} \int_{Q\cap \{U_1 >  U_2 \}} ( U_1(x) - u_{Q^{+}})^{+} \dx \dy \\
& \qquad + \dashint_{Q^{++}} (u_{Q^{+}} - U(y) )^{+} \dx \dy 
\end{align*}
where the second term is bounded by the very definition of $\PBMO^{-}$ together with the fact $u \leq U$, and the first term is dealt with the small rectangle inequality \eqref{eq:small_rectangle_final} that we proved previously. This completes the proof that $M_*^{\gamma-} u  \in  - \BMO^{-}$, and consequently we have that $M_*^{\gamma-}u \in \PBMO^{-}$.
\end{proof}

\subsection{Concluding remarks}
The preceding lemma contains essentially all that we wanted to prove. We rephrase the results in the next theorem. Note that for time-independent functions this theorem gives the classical $\BMO \to \BMO$ result whereas the the case with non-trivial time-dependency is new and interesting.
 
\begin{theorem}
\label{thm:main}
Let $u = u^{+} - u^{-} \in \PBMO^{-}$ have a locally integrable $M_*^{\gamma-}$ maximal function. Then
\begin{itemize}
\item[(i)] $M_*^{\gamma-} u^{+} \in \PBMO^{-}$
\item[(ii)] $M_*^{\gamma-}  (-u^{-}) \in \PBMO^{+}$
\item[(iii)] The maximal functions above, $U^{+}$ and $U^{-}$, control pointwise $M_*^{\gamma-} u$ in the following way:
\[\max\{U^{-},U^{+}\} \leq M_*^{\gamma-}u \leq U^{-} + U^{+} .\] 
\end{itemize}
\end{theorem}
\begin{proof}
The first item is Lemma \ref{lemma:boundedness}. The second one follows from the fact that $u^{-} \in \PBMO^{+}$ by Proposition \ref{prop:katkot}, and 
\[M_*^{\gamma-}(-u^{-}) = M_*^{\gamma +}(u^{-}) \in \PBMO^{+} \]
by symmetry. The third item is obvious by the previous ones.
\end{proof}
Functions of the type $e^{t}-e^{-t}$ on $\mathbb{R}$ show that the third item is almost the best one may hope for. The method of the proof of Lemma \ref{lemma:boundedness} does not seem to give a better result than the one above. However, already the maximal function belonging to the sum space $\PBMO^{-} + \PBMO^{+}$ would look very much like the classical $\BMO \to \BMO$ result for the Hardy-Littlewood maximal function. On the other hand, a function $u \in \BMO(\mathbb{R}^{n})$ satisfies $u \in \PBMO^{+} \cap \PBMO^{-}$ and $M_{*}^{\pm}u = M$ in the time-independent case so the well-known stationary result is covered by our evolutionary theorem. 

\section{A one-dimensional result}
The question about boundedness of the one-sided maximal operators on $\BMO^{\pm}$ spaces of Mart\'in-Reyes and de la Torre \cite{MRT1994} has also not been studied prior to this work, at least to our best knowledge. In this setting, a statement corresponding to Lemma \ref{lemma:boundedness} is much easier to prove, but we point out that the correct claim cannot be deduced directly from the multidimensional result. 

In dimension one, the relevant maximal function of $u \geq 0$ is usually defined as
\[U(x) = \sup_{h > 0} \frac{1}{h} \int_{x-h}^{x} u,\]
and there is no gap between the evaluation point $x$ and the domain of integration $(x-h,x)$. It is a general phenomenon that expressions that have a gap in the parabolic context seldom have it in dimension one. It is also usual that the gap is qualitatively inessential in dimension one whereas it is only quantitatively inessential in the parabolic context. We refer to \cite{KS2016_1} and \cite{MRT1994} for more precise discussion on what is the role of the gap in the theory of one-sided maximal functions, weights, and BMO.

As it were, a direct application of Lemma \ref{lemma:boundedness} does not give the optimal result in the one-dimensional context of \cite{MRT1994}. In order to ease the task of the reader only interested in the one-dimensional case, we give a proof of the correct one-dimensional version of Lemma \ref{lemma:boundedness}. The exposition of this proof is intended to be independent of all the other sections of this paper. 

We start by recalling the definitions from \cite{MRT1994}. We say that $u \in \BMO^{+}(\mathbb{R})$ if 
\[\norm{u}_{*} := \sup_{I} \frac{1}{\abs{I}} \int_{I} (u-u_{I^{+}})^{+}   < \infty .\]
The supremum is over all intervals, and $I^{+} = I+ \abs{I}$. This number is comparable to 
\[\sup_{I} \frac{1}{\abs{I}^{2}} \int_{I} \int_{I^{+}} (u(t_1) - u(t_2))^{+} \dt_2 \dt_1\]
as one can deduce from the theorems 2 and 3 in \cite{MRT1994}. Finally, remember the definition of the one-sided maximal function given in the beginning of this section.

\begin{theorem}
Let $u \in \BMO^{+}(\mathbb{R})$ be positive. If $U$ is locally integrable, then
\[\norm{U}_{*} \leq C \norm{u}_{*}\]
\end{theorem}
\begin{proof}
Take $u \geq 0$ from $\BMO^{+}(\mathbb{R})$. Let $I$ be an interval. We note that $U = \max(U_1, U_2 )$ where 
\begin{align*}
U_1(x) &= \sup \{ u_{(x-h,x)} : h \leq \abs{I} \} \\
U_2(x) &= \sup \{ u_{(x-h,x)} : h \geq \abs{I} \}.
\end{align*}

\subsection*{Large $h$}
We start by estimating $(U_2(x) - U_{I^{+}})^{+}$ where $I^{+} = I + \abs{I}$, and $x \in I$. Take any $h \geq \abs{I}$. Choose any $y \in I^{+}$. 
Suppose that $x \leq y - h$. Then
\begin{align*}
(u_{(x-h,x)} - u_{(y-h,y)}  )^{+} &\leq \frac{1}{h^{2}} \int_{x-h}^{\frac{x+y-h}{2}} \int_{\frac{x+y-h}{2}}^{y} (u(t_1) - u(t_2))^{+} \dt_2 \dt_1  \\
& \leq \frac{C}{h^{2}} \cdot \left( \frac{h+y-x}{2} \right)^{2} \norm{u}_{*} \leq C \norm{u}_{*}
\end{align*}
since $y-x \leq 2 \abs{I} \leq 2h$. 

In the complementary situation $x > y-h$, let $k$ be the positive integer such that $x - k (y-x) \in [x-h,y-h]$. The proof of the claim proceeds through an iterative algorithm which either stops after hitting the case $k = 1$ or converges asymptotically.

\subsection*{Case $k=1$} 
The property $k  = 1$ is equivalent to 
\[ y-x \geq \frac{1}{2} h .\]
If this holds, we can bisect both $(x-h,x)$ and $(y-h,y)$ into two halves of equal length. Call them $I_x^{\pm}$ and $I_y^{\pm}$. Then
\begin{align*}
(u_{(x-h,x)} - u_{(y-h,y)}  )^{+} &= 2 (u_{I_x^{-}} + u_{I_x^{+}} - (u_{I_y^{-}} + u_{I_y^{+}})  ) \\
& \leq 2 (u_{I_x^{-}} - u_{I_y^{-}}) + 2(u_{I_x^{+}} - u_{I_y^{+}}) \\
& \leq C \norm{u}_{*}
\end{align*}
by a reduction to the case already handled. 

\subsection*{Case $k>1$}
Now the intersection of the intervals is big. Denote $y-x = d$ and $I_j = (x-jd,x-(j-1)d)$.
Then
\begin{align*}
&(u_{(x-h,x)} - u_{(y-h,y)}  )^{+} \\
&= \frac{1}{h} \left(   \sum_{j=1}^{k} \int_{I_j}u -  \sum_{j=1}^{k} \int_{I_{j-1}} u  + \int_{x-h}^{x- kd}u  - \int_{y-h}^{x- (k-1)d}u  \right)^{+}.
\end{align*}
We modify the left-most intervals in the sums by setting
\begin{align*}
\tilde{I}_k &= I_k  \cup (x-h,x-kd) \quad \textrm{and} \\
\tilde{I}_{k-1} &= I_{k-1} \cup (y-h,x- (k-1)d).
\end{align*}
The other intervals we keep as they are. With this notation, we can continue to estimate
\begin{align*}
(u_{(x-h,x)} - u_{(y-h,y)}  )^{+} &\leq  \frac{\abs{\tilde{I}_k}}{h} (u_{\tilde{I}_k} - u_{\tilde{I}_{k-1}})^{+}  +  \frac{1}{h} \sum_{j=1}^{k-1} \abs{I_j} (u_{I_j} - u_{I_{j-1}})^{+}  \\
&\leq \frac{1}{h} \left( \abs{\tilde{I}_k}(u_{\tilde{I}_k} - u_{\tilde{I}_{k-1}})^{+} + (k-1) d  \norm{u}_{*} \right)   \\
&=  \frac{1}{h} \left( \abs{\tilde{I}_k}(u_{\tilde{I}_k} - u_{\tilde{I}_{k-1}})^{+} + (h- \abs{\tilde{I}_k}) \norm{u}_{*} \right) .
\end{align*}
Here $\tilde{I}_k \cap \tilde{I}_{k-1} \neq \varnothing$ and we can keep iterating the process.
Namely, we can repeat the argument for $(u_{\tilde{I}_k} - u_{\tilde{I}_{k-1}})^{+}$. If $|\tilde{I}_k \cap \tilde{I}_{k-1}| \leq \frac{1}{2} | \tilde{I}_k |$ the process terminates after an application of Case $k=1$ on these intervals. Otherwise we iterate the process as follows.

Denote
\begin{align*}
(x-h,x) &= I_{l}^{(1)} \quad \textrm{and} \quad (y-h,y) = I_{r}^{(1)} \\
\tilde{I}_k &= I_{l}^{(2)} \quad \textrm{and} \quad \tilde{I}_k = I_{r}^{(2)} 
\end{align*}
and form recursively the intervals $(I_l^{(i)}, I_r^{(i)})$ for all $i \in \{2,\ldots, K\}$ where $K$ is the value of $i$ at which
\[ |I_l^{(i)} \cap I_r^{(i)}| \geq \frac{1}{2} |I_l^{(i)} |\]
is violated. At that point we are done by the case $k = 1$. It is also possible that $K = \infty$. Then the process does not stop, and we have that 
\[ \abs{I}_{l}^{(i)} (u_{I_l^{(i)}} - u_{I_r^{(i)}})^{+} \longrightarrow 0 \]
as $i \to \infty$ by local integrability of $u$ since $\abs{I_l^{(i)}} < |I_l^{(i-1)}|$ for all $i$. There are also numbers $\theta_i = \abs{I_l^{(i)}} - \abs{I_l^{(i+1)}}  \geq 0$ such that
\[ \sum_{i=1}^{K} \theta_i \leq h \]
and
\begin{align*}
h(u_{(x-h,x)} - u_{(y-h,y)}  )^{+} &\leq \abs{\tilde{I}_k}(u_{\tilde{I}_k} - u_{\tilde{I}_{k-1}})^{+} + (h- \abs{\tilde{I}_k}) \norm{u}_{*} \\
& = \abs{\tilde{I}_k}(u_{\tilde{I}_k} - u_{\tilde{I}_{k-1}})^{+} + \theta_1 \norm{u}_{*} \\
& \leq C \norm{u}_{*} h + \norm{u}_{*} \sum_{i=1}^{\infty} \theta_i \\
& \leq C h \norm{u}_{*}.
\end{align*}

Consequently, we always have 
\[(u_{(x-h,x)} - U(y)  )^{+} \leq C \norm{u}_{*}\]
for all $y \in I^{+}$. Taking the supremum over $h \geq \abs{I}$, we see that
\[\frac{1}{\abs{I}} \int_{I} (U_2 - U_{I^{+}}) \leq C \norm{u}_{*}. \]

\subsection*{Small $h$}
Then we move to the part dealing with $U_1$. Now the maximal function averages over small intervals. To control the quantity
\[\int_{I} (U_{1} - u_{I^{+} \cup I^{2+}})^{+} \]
where $I^{2+}=I + 2 \abs{I}$, we form a ``one-sided'' Calder\'on-Zygmund decomposition of $u$ in $I \cup I^{-}$ where $I^{-} = I - \abs{I}$ . More precisely, let $\lambda > u_{I^{+}\cup I^{2+}}$. Take the maximal dyadic subintervals $\{I_i\}_i$ of $I \cup I^{-}$ with $u_{I_i^{+}}  > \lambda$. Let $\widehat{I_i}$ be the parent of $I_i$. We write
\begin{align*}
b_i &= 1_{I_i} (u-u_{\widehat{I_i}^{+}}) \\
b &= \sum_{i} b_i \\
g &= \sum_i 1_{I_i} u_{\widehat{I_i}^{+}} + 1_{I \setminus \cup_i I_i} u.
\end{align*}
Here $g \leq \lambda$. By the definition of $U_1$, we have that everywhere in $I$ it holds that
\[U_1 \leq M(1_{I^{-} \cup I}u) \leq M(b^{+}) + M g.   \]
Here $M$ is the standard two-sided Hardy-Littlewood maximal operator. Since $b_i$ have disjoint supports, we can use the John-Nirenberg inequality of \cite{MRT1994} to estimate 
\begin{align*}
\int_{\mathbb{R}} (b^{+})^{2} &=  \sum_i \int_{I_i} (u-u_{\widehat{I_i}^{+}})_{+}^{2} \leq \sum_i \abs{\widehat{I_i}}  \frac{1}{\abs{\widehat{I_i}}}\int_{\widehat{I_i}} (u-u_{\widehat{I_i}^{+}})_{+}^{2} \leq C \abs{I} \norm{u}_{*}^{2}.
\end{align*}
Denote $\Omega = I \cap \{U_{1} > u_{I^{+} \cup I^{2+}} \}$. Then
\begin{align*}
\int_{\Omega} U_{1} &\leq \abs{\Omega}^{1/2} \norm{M(b^{+})}_{L^{2}} + \abs{\Omega} \lambda \\
& \leq C \abs{I} \norm{u}_{*} + \abs{\Omega} \lambda,
\end{align*}
which proves that 
\[\int_{I} (U_{1} - u_{I^{+} \cup I^{2+}})^{+} \leq C \norm{u}_{+} \abs{I}.\]

Finally, for any interval $J$ we may write $J = I^{-} \cup I$ where $I^{-} = I - \abs{I}$. Then
\begin{align*}
&\int_{J} (U_{1} - u_{ J^{+}})^{+}  = \int_{I^{-}} (U_{1} - u_{I^{+} \cup I^{2+}})^{+} + \int_{I} (U_{1} - u_{I^{+} \cup I^{2+}})^{+} \\
&= I + II
\end{align*}
where $II \leq C \norm{u}_{*} \abs{I}$ by the previous considerations. For the other part, we estimate
\begin{align*}
\int_{I^{-}} (U_{1} - u_{I^{+} \cup I^{2+}})^{+} & \leq \int_{I^{-}} (U_{1} - u_{I \cup I^{+}})^{+} + \abs{I} ( u_{I \cup I^{+}} - u _{I^{+} \cup I^{2+}})^{+} \\
& \leq C \norm{u}_{*} \abs{I} + \frac{\abs{I}}{2}(u_{I}+u_{I^{+}} - u _{I^{+}} - u_{ I^{2+}}  )^{+}  \\
& \leq C \norm{u}_{*} \abs{I} + \frac{\abs{I}}{2} [(u_{I} - u _{I^{+}})^{+}  + (u_{I^{+}}  - u_{ I^{2+}}  )^{+}] \\
& \leq (C+1) \norm{u}_{*} \abs{I}.
\end{align*}
This proves the estimate for the part of $U_1$.

Putting the pieces together, we see that 
\[\int_{I} (U  - U_{I^{+}})^{+} \leq C \norm{u}_{*} \abs{I}\]
for all intervals $I$ and a numerical $C$. The proof is complete.
\end{proof}



\end{document}